\documentclass[12pt]{amsart}

\usepackage{graphicx}
\usepackage{amsmath}
\usepackage{amssymb}
\usepackage{amsfonts}
\usepackage{verbatim}
\usepackage{scalefnt}
\usepackage{multirow}
\usepackage{enumerate}
\usepackage{mathtools}
\usepackage{float}
\usepackage{enumitem}
\restylefloat{figure}
\usepackage[margin=3cm]{geometry}

\usepackage{fancyhdr}
\usepackage{hyperref}

\newtheorem{theorem}{Theorem}

\newtheorem{corollary}[theorem]{Corollary}

\newtheorem*{definition*}{Definition}

\newtheorem{remark}[theorem]{Remark}

\numberwithin{equation}{section}
\numberwithin{theorem}{section}

\title[Product representations of perfect powers]%
  {Product representations of perfect powers}

\author{P\'eter P\'al Pach}
\email{pachpp@renyi.hu}
\address{HUN-REN Alfr\'ed R\'enyi Institute of Mathematics, Re\'altanoda utca 13--15., H-1053 Budapest,  Hungary; \newline \hspace*{4mm}
MTA--HUN-REN RI Lend\"ulet ``Momentum'' Arithmetic Combinatorics Research Group, Re\'altanoda utca 13--15., H-1053 Budapest,  Hungary; \newline \hspace*{4mm}
Department of Computer Science and Information Theory, Budapest University of Technology and Economics, M\H{u}egyetem rkp. 3., H-1111 Budapest, Hungary;
\newline \hspace*{4mm} 
Extremal Combinatorics and Probability Group (ECOPRO), Institute for Basic Science (IBS), Daejeon, South Korea.}

\author{Csaba S\'andor}
\email{sandor.csaba@ttk.bme.hu}
\address{Department of Stochastics, Institute of Mathemetics, Budapest University of Technology and Economics, M\H{u}egyetem rkp. 3., H-1111 Budapest, Hungary; \newline \hspace*{4mm}
HUN-REN Alfr\'ed R\'enyi Institute of Mathematics, Re\'altanoda utca 13--15., H-1053 Budapest,  Hungary; \newline \hspace*{4mm}
MTA--HUN-REN RI Lend\"ulet ``Momentum'' Arithmetic Combinatorics Research Group, Re\'altanoda utca 13--15., H-1053 Budapest,  Hungary.}

\thanks{}

\begin{document}

\begin{abstract}	
Let $\rho_k(N)$ denote the maximum size of a set $A\subseteq \{1,2,\dots,N\}$ such that no product of $k$ distinct elements of $A$ is a perfect $d$-th power. In this short note, we prove that $\rho _d(N)=\sum\limits_{k=1}^{d-1}\pi\left( \frac{N}{k} \right) +O_d(\pi (N^{1/2}))$, furthermore, for prime power $d$ and sufficiently large $N$ we have $\rho _d(N)=\sum\limits_{k=1}^{d-1}\pi\left( \frac{N}{k} \right)$. This answers a question of Verstra\"ete.
\end{abstract}

\date{\today}
\maketitle

\section{Introduction}

Let $\rho _d(N)$ denote the maximum size of a subset of  $[N]:=\{1,2,\dots,N\}$ with no product representation of the polynomial $f(x)=x^d$, that is $a_{1}a_{2}\dots a_{k}\ne x^d$ for any distinct $a_{1},a_2,\dots,a_k\in A$ and integer $x$. It is easy to see that $\rho _2 (N)=\pi (N)$, the number of primes up to $N$, for every positive integer $N$. (Indeed, the exponent vectors in the canonical forms of elements for such sets form a linearly independent system over the 2-element field $\mathbb{F}_2$).

In 2006, Verstra\"ete~\cite{Ver} proved the following bound.
\begin{theorem}[Verstra\"ete]
Let $d$ be a prime number. Then
$$
\rho _d(N)=\sum_{k=1}^{d-1}\pi\left( \frac{N}{k} \right) +O(\pi (N^{1/2})).
$$
\end{theorem}
He wrote that it would be interesting to determine the order of magnitude of $$
\rho _d(N)-\sum_{k=1}^{d-1}\pi\left( \frac{N}{k} \right)
 $$
 in the above theorem, and also mentioned that perhaps it is always at most a constant.

We answer this question in the following stronger form.

\begin{theorem}\label{thm-exact}
Let $d\ge 2$ be a prime power. For $N\ge 2^dp_d$, $$\rho _d(N)= \sum_{k=1}^{d-1}\pi \left( \frac{N}{k}\right),$$
where $p_d$ denotes the $d$-th prime.
\end{theorem}

\begin{remark}
For very small values of $N$ the above equality may not hold. 

A calculation shows that for $d=3$, we have $\rho _3(N)= \sum\limits_{k=1}^2\pi \left( \frac{N}{k}\right)$ if $N\ge 15$, but $9=\rho _3(14)< \sum\limits_{k=1}^2\pi \left( \frac{14}{k}\right)=10,$ therefore the strongest statement for $d=3$ is that $\rho _3(N)= \sum\limits_{k=1}^2\pi \left( \frac{N}{k}\right)$ for $N\ge 15$.

Clearly,
$$
\sum\limits_{k=1}^{d-1}\pi \left( \frac{2d-2}{k}\right)=\sum\limits_{k=1}^{d-1} \sum_{p\le \frac{2d-2}{k}}1= \sum_{p\le 2d-2}\sum_{k\le \frac{2d-2}{p}}1=\sum_{p\le 2d-2}\left\lfloor \frac{2d-2}{p}\right\rfloor =\sum_{n\le 2d-2}\omega (n),
$$
where $\omega (n)$ is the number of distinct prime divisors of the positive integer $n$. For $d\ge 4$, we can prove by induction on $d$ that $\sum\limits_{k=1}^{d-1}\pi \left( \frac{2d-2}{k}\right)\ge 2d-2$. 

If $A\subseteq {N}$ has no product representation of the polynomial $f(x)=x^d$, then $1\not \in A$. It follows that $\rho _d(N)<N$ for every positive integer $N$. Hence, $\rho _d(2d-2)< \sum\limits_{k=1}^{d-1}\pi \left( \frac{2d-2}{k}\right)$ for every $d\geq 4$.

So for every $d\ge 3$ there exists an integer $N_d$ such that $\rho _d(N_d)< \sum\limits_{k=1}^{d-1}\pi \left( \frac{N_d}{k}\right)$.
\end{remark}

For general integer $d$, we could prove a bound with the same error term that Verstra\" ete had for primes. 

\begin{theorem}\label{thm-upp-gen}
Let $d\ge 2$ be an integer. Then
$$
\rho _d(N)=\sum_{k=1}^{d-1}\pi\left( \frac{N}{k} \right) +O_d(\pi (N^{1/2})).
$$
\end{theorem}

We shall use the Davenport constant in our proofs. 

Let $G$ be a finite Abelian group. The Davenport constant $D(G)$ is defined as the smallest positive integer $\ell$ such that every sequence of $\ell$ elements contains a non-empty subsequence adding up to 0.

The exact value of the Davenport constant is known only for a few classes of Abelian groups, but we know it precisely for $p$-groups. Let $p$ be a prime number. In 1969, Olson \cite{O} proved that if $G=\oplus _{i=1}^r\mathbb{Z}_{p^{e_i}}$, then $$\displaystyle D(G)=1+\sum_{i=1}^r(p^{e_i}-1).$$ 

Let $\exp (G)$ be the exponent of the Abelian group $G$. For any finite Abelian group we have (see \cite{AGP,EK,M})
$$
D(G)\le \exp(G)\left( 1+\log \frac{|G|}{\exp (G) } \right).
\eqno(*)$$

For general integer $d$, we prove the following estimates.

\begin{theorem}\label{upper1}
Let $d\ge 2$. For every positive integer $N$,

$$\rho _d(N)\le \left[\sum_{k=1}^{d-1}\pi \left( \frac{N}{k}\right) \right] +D(\mathbb{Z}_d^{\pi \left( {N}/{d} \right) })- (d-1)\pi \left( {N}/{d}\right) -1.$$
\end{theorem}

\begin{theorem}\label{upper2}
Let $d\ge 2$. For every positive integer $N$,
$$\rho _d(N)\le \left[\sum\limits_{k=1}^{d-1}\pi \left( \frac{N}{k}\right)\right] +dD\left(\mathbb{Z}_d^{\pi \left( \sqrt{N} \right) }\right).$$
\end{theorem}
Using the bound $(*)$ Theorem~\ref{upper2} yields the following corollary verifying Theorem~\ref{thm-upp-gen}.
\begin{corollary}
For $d\ge 2$, we have $$ \rho_d (N)\le \sum_{k=1}^{d-1}\pi \left( \frac{N}{k}\right) +d\log d \cdot \pi (\sqrt{N}).$$
\end{corollary}

%%%% Notations
\textbf{Notations.} Throughout the paper, we denote by $[N]$ the set $\{1,2, \ldots, N\}$. The standard notation $O$ is applied to positive quantities in the usual way.
That is, $Y = O(X)$ means that $X \geq cY$,
for some absolute constant $c > 0$. 
If the constant $c$ depends on a quantity
$t$, we write $Y = O_t(Y )$. 

The primes are denoted by $p_1=2,p_2=3,\dots$ in the paper.

\section{Proofs}

\begin{proof}[Proof of Theorem~\ref{upper1}]
Let $A\subseteq [N]$ be a set avoiding product representations of $x^d$. For $1\le k\le d-1$, define 
$$A_k=\left\{ a\in A: \exists \text{ prime } p: p|a,\  \frac{N}{k+1}< p\leq\frac{N}{k}  \right\},$$
notice that $|A_k|\le k\left( \pi \left( \frac{N}{k}\right) -\pi \left( \frac{N}{k+1} \right) \right) $, since each prime $p\in \left(\frac{N}{k+1},\frac{N}{k}\right]$ has only $k$ multiples in $[N]$.

Finally, let $A_d$ be the  set of those $a\in A$ whose prime divisors are all at most $\frac{N}{d}$. Let the canonical form of $a\in A_d$ (allowing 0 exponents) be
$$a=p_1^{\alpha _1}p_2^{\alpha _2}\dots p_{\pi ({N}/{d})}^{\alpha _{\pi ({N}/{d})}},$$
let us assign to $a$ the vector 
$$( \alpha_1 \pmod{d},\alpha _2 \pmod{d},\dots ,\alpha _{\pi ({N}/{d})}\pmod{d} )\in \mathbb{Z}_d^{\pi ({N}/{d})}.$$
Since these vectors form a multiset that avoids zero sums in $\mathbb{Z}_d^{\pi ({N}/{d})}$,  we have $|A_d|\le D(\mathbb{Z}_d^{\pi ({N}/{d})})-1$. Hence,
\begin{multline*}
|A|\le \sum_{k=1}^d|A_k|\le \sum_{k=1}^{d-1}k\left(\pi \left( \frac{N}{k}\right) -\pi \left( \frac{N}{k+1} \right) \right) + D(\mathbb{Z}_d^{\pi ({N}/{d})})-1 =
\\
=\sum_{k=1}^{d-1}\pi \left( \frac{N}{k} \right) + D(\mathbb{Z}_d^{\pi ({N}/{d})})-1 -(d-1)\pi \left( {N}/{d} \right).
\end{multline*}
\end{proof}

\begin{proof}[Proof of Theorem~\ref{upper2}]
Let $A\subseteq [N]$ be a set avoiding product representations of $x^d$. Let us partition the elements of $A$ into three parts: $A=Q\cup R\cup S$, where
\begin{itemize}
\item each element of $Q\cup R$ has a prime divisor larger than $\sqrt{N}$,
\item each prime divisor of elements from $S$ is at most $\sqrt{N}$,
\item for every prime $p>\sqrt{N}$ the number of multiples of $p$ contained in $R$ is divisible by $d$,
\item for every prime $p>\sqrt{N}$ the number of multiples of $p$ contained in $Q$ is at most  $d-1$.
\end{itemize}
That is, we put the $\sqrt{N}$-smooth elements of $A$ to $S$, and for each prime $p>\sqrt{N}$ if the number of multiples of $p$ in $A$ is $k_p$, then we put $k_p\mod {d}$ multiples of $p$ to $Q$ and the rest of the multiples to $R$. (Note that every element of $A\subseteq [N]$ can have at most one such prime divisor.)

A calculation similar to the one in the proof of Theorem~\ref{upper1} shows that 
$$|Q|\leq \left[\sum\limits_{k=1}^{d-1} \pi\left(\frac{N}{k}\right)\right]-(d-1)\pi(\sqrt{N}),$$
giving the main term of the estimate.

To bound the size of the set $R\cup S$, for each prime $p>\sqrt{N}$ let us partition the multiples of $p$ in the set $R$ into groups of size $d$. In such a group the product of the $d$ elements is of the form $p^dt$, where $t$ is the $\sqrt{N}$-smooth part: $t$ only has prime divisors being at most $\sqrt{N}$. Let $T$ be the multiset containing these $\sqrt{N}$-smooth parts and the elements of $S$. Note that elements of the multiset $T$ only has prime divisors being at most $\sqrt{N}$. Let the canonical form of $t\in T$ (allowing 0 exponents) be
$$t=p_1^{\alpha _1}p_2^{\alpha _2}\dots p_{\pi (\sqrt N)}^{\alpha _{\pi (\sqrt N)}},$$
let us assign to $t$ the vector 
$$( \alpha_1 \pmod{d},\alpha _2 \pmod{d},\dots ,\alpha _{\pi (\sqrt N)}\pmod{d} )\in \mathbb{Z}_d^{\pi (\sqrt N)}.$$
Since these vectors form a multiset that avoids zero sums in $\mathbb{Z}_d^{\pi (\sqrt N)}$,  we have 
$$|T|\leq D(\mathbb{Z}_d^{\pi (\sqrt N)})-1.$$
Hence,
$$|A|=|Q|+|R|+|S|\leq |Q|+d|T|\leq \left[\sum\limits_{k=1}^{d-1}\pi \left( \frac{N}{k}\right)\right]+ dD(\mathbb{Z}_d^{\pi(\sqrt N)}). $$

\end{proof}

\begin{proof}[Proof of Theorem~\ref{thm-exact}]
First of all, since, for prime power $d$, we have 
$$D\left(\mathbb{Z}_d^{\pi \left( {N}/{d} \right) }\right)=(d-1)\pi \left( {N}/{d}\right)+1,$$ 
from Theorem~\ref{upper1} we immediately get the required upper bound
$$
\rho _d(N)\le \sum_{k=1}^{d-1}\pi \left( \frac{N}{k}\right).
$$
To finish the proof we present a matching lower bound construction.

For $1\le k\le d-2$, let 
$$A_k=\left\{ip:  \frac{N}{k+1}<p\le \frac{N}{k},\ p\text{ prime},\ 1\le i\le k\right\},$$ 
and let 
$$A_{d-1}=\left\{ip: p_d< p\le \frac{N}{d-1},\ p\text { prime},\quad 1\le i\le d-1\right\}.$$ 
Choose the integers $0\le j_1<j_2<\dots <j_{d}\le d$ such that $j_1+\dots +j_d \equiv 1 \pmod{d}$. 

For $2\le  s\le d$, let us set
$$
B_{p_s}=\{ 2^{j_t}p_s: 1\le t\le d\} .
$$
Finally, let
$$
A:=\bigcup\limits_{1\le k\le d-1}A_k \cup \left( \bigcup\limits_{s=2}^d B_{p_s}\right) .
$$
Now, we show that there is no product representation of a perfect $d$-th power in the set $A$. For the sake of contradiction, let us assume that $a_1a_2\dots a_{\ell}=x^d$ for some $a_i\in A$, $a_1<a_2<\dots <a_{\ell}$ and $x\in \mathbb{Z}^+$. 

If $a_j\in A_k$ for some $1\le k\le d-1$, then $a_j=ip$, where $1\le i\le d-1$, then the exponent of $p$ in the product $a_1a_2\dots a_{\ell}$ is between $1$ and $d-1$, a contradiction. 

Thus we may assume that $a_j\in \bigcup\limits_{s=2}^d B_{p_s}$ for every $1\le j\le \ell$. If $a_j\in B_{p_s}$, $2\le s\le d$, then the exponent of $p_s$ in the product $a_1a_2\dots a_{\ell}$ is between $1$ and $d$, so it has to be exactly $d$, implying that all the $d$ multiples of $p_d$ were chosen in the set $\{a_1,\dots,a_{\ell}\}$. 

If this happens for $u$ such primes, say, for $p_{s_1},\dots,p_{s_u}$ where $2\le s_1<s_2<\dots <s_u\le d$, then $a_1a_2\dots a_{\ell}=2^vp_{s_1}^d\dots p_{s_u}^d$, and $v\equiv u \pmod{d}$ according to the definition of the sets $B_{p_s}$. However, $1\leq u\leq d-1$, which is a contradiction.

%If $a_1a_2\dots a_l=2^vp_{s_1}^d\dots p_{s_u}^d$, $2\le s_1<s_2<\dots s_u\le d$, then $v\equiv u \pmod{d}$, where $1\le u\le d-1$, a contradiction.

Clearly, each element of $\bigcup\limits_{s=2}^d B_{p_s}$ is at most $2^dp_d$. Hence 
\begin{eqnarray*}
|A|=&&\sum_{k=1}^{d-1}|A_k|+\sum_{2\le s\le d}|B_{p_s}|\\
=&& \sum_{k=1}^{d-2}k\left( \pi \left( \frac{n}{k}  \right)- \pi \left( \frac{n}{k+1}  \right) \right) +(d-1)\left( \pi \left( \frac{n}{d-1}  \right)- d \right) +d(d-1)\\
=&&\sum_{k=1}^{d-2}k\left( \pi \left( \frac{n}{k}  \right)- \pi \left( \frac{n}{k+1}  \right) \right) +(d-1)\pi \left( \frac{n}{d-1}  \right)\\
=&&\sum_{k=1}^{d-1}\pi \left( \frac{N}{k}\right).
\end{eqnarray*}
\end{proof}

\section{Concluding remarks and open questions}

In this note it is shown that for prime power $d\ge 2$ and $N\ge 2^dp_d$, we have $\rho _d(N)= \sum\limits_{k=1}^{d-1}\pi \left( \frac{N}{k}\right).$ It remains open whether the threshold $2^dp_d\sim 2^dd\log d$ can be improved.

For arbitrary integer $d$ we proved the bound 
$\rho _d(N)=\sum\limits_{k=1}^{d-1}\pi\left( \frac{N}{k} \right) +O_d(\pi (N^{1/2})).$ However, we believe that the error term $O_d(\pi (N^{1/2}))$ can be further improved.

\section{Acknowledgements}
Both authors were supported by the National Research, Development and Innovation Office NKFIH (Grant Nr. K146387). PPP was also supported by the National Research, Development and Innovation Office NKFIH (Excellence program, Grant Nr. 153829) and by the Institute for Basic Science (IBS-R029-C4).

\medskip

\end{document}